\documentstyle[12pt]{article}
\newtheorem{theorem}{Theorem}[section]

\newtheorem{corollary}[theorem]{Corollary}
\newtheorem{proposition}[theorem]{Proposition}
\newtheorem{remark}[theorem]{Remark}
\newtheorem{definition}[theorem]{Definition}
\newtheorem{example}[theorem]{Example}
\def\qed{\hfill{$\Box$}}
\def\proof{\indent{\bf Proof. }}
\def\kappa{K_{S}^{2}}

\def\bbbc{{\mathchoice {\setbox0=\hbox{$\displaystyle\rm C$}\hbox{\hbox
 to0pt{\kern0.4\wd0\vrule height0.9\ht0\hss}\box0}}
 {\setbox0=\hbox{$\textstyle\rm C$}\hbox{\hbox
 to0pt{\kern0.4\wd0\vrule height0.9\ht0\hss}\box0}}
 {\setbox0=\hbox{$\scriptstyle\rm C$}\hbox{\hbox
 to0pt{\kern0.4\wd0\vrule height0.9\ht0\hss}\box0}}
 {\setbox0=\hbox{$\scriptscriptstyle\rm C$}\hbox{\hbox
 to0pt{\kern0.4\wd0\vrule height0.9\ht0\hss}\box0}}}}

 \def\bbbz{{\mathchoice {\hbox{$\sf\textstyle Z\kern-0.4em Z$}}
 {\hbox{$\sf\textstyle Z\kern-0.4em Z$}}
 {\hbox{$\sf\scriptstyle Z\kern-0.3em Z$}}
 {\hbox{$\sf\scriptscriptstyle Z\kern-0.2em Z$}}}}

 \def\bbbq{{\mathchoice {\setbox0=\hbox{$\displaystyle\rm Q$}\hbox{\raise
 0.15\ht0\hbox to0pt{\kern0.4\wd0\vrule height0.8\ht0\hss}\box0}}
 {\setbox0=\hbox{$\textstyle\rm Q$}\hbox{\raise
 0.15\ht0\hbox to0pt{\kern0.4\wd0\vrule height0.8\ht0\hss}\box0}}
 {\setbox0=\hbox{$\scriptstyle\rm Q$}\hbox{\raise
 0.15\ht0\hbox to0pt{\kern0.4\wd0\vrule height0.7\ht0\hss}\box0}}
 {\setbox0=\hbox{$\scriptscriptstyle\rm Q$}\hbox{\raise
 0.15\ht0\hbox to0pt{\kern0.4\wd0\vrule height0.7\ht0\hss}\box0}}}}
 
\input xy
\xyoption{all}

 1

\font\es=msbm10 scaled\magstep 1

\def\Natural{\mbox{\es N}}

\def\Rational{\mbox{\es Q}}

\def\proj{\mbox{\es P}}

\begin{document}

\vglue2truecm

\begin{center}
{\large\bf ON THE SLOPE OF FIBRED SURFACES}\\
\end{center}

\vglue1truecm

\begin{center}
Miguel Angel BARJA\footnotemark[1]\footnotetext[1]
{Partially supported by
CICYT PS93-0790 and HCM project n.ERBCHRXCT-940557.}\\
{\scriptsize{Departament de Matem\`atica Aplicada I\\
Universitat Polit\`ecnica de Catalunya\\
Diagonal 647\\
08028 Barcelona. Spain\\
e-mail: barja@ma1.upc.es}}

\smallskip

Francesco ZUCCONI\footnotemark[2]\footnotetext[2]
{Partially supported by HCM project n.ERBCHRXCT-940557.}\\
{\scriptsize{Dipartamento di Matematica e Informatica\\
Universit\`a degli studi di Udine\\
Via delle Scienze, 206\\
33100 Udine. Italy\\
e-mail: zucconi@dimi.uniud.it}}

\end{center}

\setcounter{section}{-1}

\section{Introduction}

Let $f: S \longrightarrow B$ be a projective, surjective morphism
from a complex smooth surface $S$ onto a complex smooth curve $B$. We set $F$ for
the general fibre of $f$ and assume it is connected. Let $g=g(F)$
and $b=g(B)$. We will assume that $f$ is relatively minimal, i.e.,
that there is no $(-1)-$rational curve on fibres. We usually call $f$ a
{\it fibration} or a {\it minimal genus $b$ pencil of curves of genus $g$}.
We say that $f$ is {\it smooth} if all its fibres are smooth,
that it is {\it isotrivial} it all its smooth fibres are reciprocally
isomorphic, and that it is {\it locally trivial} if it is
smooth and isotrivial.

Our results enable to study the geographical problem of $f$
(that is, to relate numerical invariants of $F$, $S$ and $B$)
through the control of some geometrical properties of the general fibre $F$
or the influence of some {\it global} properties of $S$ such as the
relative irregularity $q_{f}=q-b$. Now we recall the basic relative
invariants for $f$. Let $\omega_S$, $\omega_B$ be the canonical
sheaves and $K_S$, $K_B$ canonical divisors of $S$ and $B$ respectively.
As usual we consider $p_g=h^0(S,\omega_S)$, $q=h^1(S,\omega_S)$,
$\chi {\cal O}_S=p_g-q+1$ and denote as $e(X)$ the topological
Euler characteristic of X. Then we set:
\begin{eqnarray}
&&K^{2}_{S/B}=(K_{S}-f^{\ast}K_{B})^{2}=K^{2}_{S}-8(b-1)(g-1)\nonumber\\
&&\chi_{f}=\mbox {deg}f_{\ast}\omega_{S/B}=
\chi {\cal O}_{S}-(b-1)(g-1)\nonumber\\
&&e_{f}=e(S)-e(B)e(F)=e(S)-4(b-1)(g-1).\nonumber
\end{eqnarray}

We have the following classical results:

\bigskip
\indent
\begin{theorem}
Let $f:S \longrightarrow B$ be a minimal genus $b$ pencil of curves
of genus $g \geq 2$. Then

\begin{itemize}

\item[(i)] (Noether) $12 \chi_f=e_f+K^2_{S/B}.$

\item[(ii)] (Zeuthen-Segre) $e_f \geq 0$. Moreover, $e_f=0$ if and only
if $f$ is smooth.

\item[(iii)] (Arakelov) $K^2_{S/B} \geq 0$. Moreover, $K^2_{S/B}=0$
if and only if $f$ is isotrivial.

\item[(iv)] $\chi_f \geq 0$. Moreover, $\chi_f=0$ if and
only if $f$ is locally trivial.

\end{itemize}

\end{theorem}

\medskip

When $f$ is not locally trivial, Xiao (cf. \cite{X}) defines the
slope of $f$ as

$$\lambda(f)=\frac{K^2_{S/B}}{\chi_f}.$$

It follows immediately from Noether's equality that $0 \leq \lambda(f)
\leq 12$. We are mostly concerned with a lower bound of the slope. The main
known result is:

\bigskip

\begin{theorem} (Cornalba-Harris, Xiao).
If $g \geq 2$ and $f$ is not locally trivial, then $\lambda(f) \geq
4-\frac{4}{g}$.

\end{theorem}

After that, the first problem was to investigate the influence of
some properties of the fibration on the behaviour of the slope.
The first direction is to study the influence of the relative
irregularity $q_f=q(S)-b$. The main known result is:

\bigskip

\begin{theorem} (Xiao \label{sabi})
If $q>b$ then $\lambda(f) \geq 4$. If $\lambda(f)=4$ and $q>b$ then
$q=b+1$ and $f_{\ast}\omega_{S/B}$ is semistable.

\end{theorem}

\medskip

The other most considered problem is the study of how properties of the
general fibre $F$ influence. Mostly due to the work of Konno
(cf. \cite{K1},\cite{K3}; see also \cite{Chen} and
\cite{SF} for other references) we know
the Clifford index (or the gonality) of the general fibre has some
meaning in the lower bound of the slope. There are evidences
for this. For example, it is known that equality $\lambda(f)=
4-\frac{4}{g}$ only holds when $F$ is hyperelliptic.
When $F$ is trigonal, a better bound is known (see Remark 3.5).
On the other hand the following theorem shows that if $g>>0$ and
$F$ has general Clifford Index then $\lambda(f)\simeq 6$ and so, together
with \ref{sabi}, it forces to understand the case $q_{f}>0$, $g>>0$
and $4<\lambda(f)<6$.

\bigskip

\begin{theorem} (Harris-Eisenbud, Konno). Assume that $g$ is odd
and that the general fibre $F$ is of general Cifford index. If
$f$ is not a semistable fibration, assume that Green's conjecture
on syzygies of the canonical curves holds. Then

$$\lambda(f) \geq 6\frac{g-1}{g+1}.$$

\end{theorem}

\medskip

An explicit sharp lower bound for $\lambda(f)$ depending on the
Clifford index of $F$ should not be easy and should depend
on other parameters. Indeed, when the fibre is trigonal, the
behaviour of $\lambda(f)$ depends on the fact that $F$ is
{\it general} or not in the set of trigonal curves (see \cite{SF}).
A similar behaviour should hold for tetragonal fibrations: there
are tetragonal fibrations with $\lambda(f)=4$ for any genus $g$
(see \cite{Ba}), although the general fibre is always bielliptic.

In \S 1 we give an idea of Xiao's method, which is our main tool.
In \S 2 we study the behaviour of $\lambda(f)$ when the general fibre
$F$ is a double cover in such a way that extends to a double cover of
the fibration itself. In \S 3 we explicit the influence of $q_{f}$ on
$\lambda(f)$ through an increasing function on
$q_f$ generalising Theorem \ref{sabi}. As natural by-product
of this estimates and the previous theorems it seems possible to construct through $\lambda(f)$
a geography for fibrations as in the case of surfaces of general type.

We set aside the double cover case by two reasons: $1)$ all double
covers are special curves of non general Clifford index;
as happens in the case of bielliptic curves with respect to tetragonal
ones, fibrations which are double
covers are candidates to be exceptional in the study of $\lambda(f)$.
$2)$ Xiao's method works very well if the possibility for $f$ to be a
double cover is excluded, which also suggest the study of double covers
as exceptional. For these fibrations we get:

\medskip

\begin{theorem} Let $f:S\rightarrow B$ be a genus $g$,
relatively minimal, non isotrivial fibration. Assume $f$ is a double cover fibration of a
fibration of genus $\gamma$. Then, if $g\geq 4\gamma+1$ we have
$$
\lambda(f)\geq 4+4\frac{(\gamma-1)(g-4\gamma-1)}{(g-4\gamma -1)(g-\gamma)+2(g-1)\gamma^{2}}.
$$
\end{theorem}

In our paper there is also a refined version of this theorem
but it is more complicate to state: see 2.4. Next theorem
is a partial answer to a natural question: how special are fibrations with
$\lambda(f)<4$?

\medskip

\begin{theorem}
Let $f:S\rightarrow B$ be a relatively minimal, non isotrivial double
cover fibration of  $\sigma:V\rightarrow B$. Let $F$ and $E$ be the fibres of $f$
and $\sigma$ respectively and let $g=g(F)$, $\gamma=g(E)$. Assume $F$ is not
hyperelliptic nor tetragonal, $\gamma\geq 1$ and $g\geq 2\gamma+11$. Then
$\lambda(f)\geq 4$.
\end{theorem}

\medskip

Our next theorem gives an affirmative answer to the expected
influence of $q_f$ on the slope as suggested by Theorem 0.3.
It gives a bound which is (strictly) increasing and
in some cases is assimptotically sharp (see Example \ref{ana}).
If $f$ is general, that is, it is not a double cover fibration we have:

\begin{theorem}
Let $f:S \longrightarrow B$ be a relatively minimal fibration which is not
a double cover fibration. Assume $g=g(F)\geq 5$ and that $f$ is not
locally trivial. Let $h=q(S)-b\geq 1$.

(i) If $h \geq 2$ and $g \geq \frac{3}{2}h+2$ then
$$\begin{array}{rl}
\quad \lambda(f)&\geq \displaystyle{
\frac{8g(g-1)(4g-3h-10)}{8g(g-1)(g-h-2)+3(h-2)(2g-1)}}
\quad \mbox{if}\ F \ \mbox{is not trigonal}\\
&\\
\quad \lambda(f)&\geq \displaystyle{
\frac{4g(g-1)(4g-3h-10)}{4g(g-1)(g-h-2)+(g-4)(2g-1)}}
\quad \mbox{if}\ F \ \mbox{is trigonal}
\end{array}$$

(ii) if $g < \frac{3}{2}h+2$ then
$$\lambda(f)\geq \frac{4g(g-1)(2g-7)}{\frac{4}{3}g(g-1)(g-3)+(g-4)(2g-1)}
$$
\end{theorem}

\medskip

As an application we obtain a nice relation between $\lambda(f)$ and the existence of other fibration on $S$ onto curves of genus at least $2$ (see also Corollaries \ref{adria1}, \ref{adria2}):

\bigskip

\begin{theorem}
Let $f:S \longrightarrow B$ be a relatively minimal, non locally trivial
fibration. Let $F$ be a fibre of $f$, $g=g(F)$ and $q=q(S)$.
Assume $f$ is not a
double cover fibration and that $s=q-b\geq1$ (i.e., $f$ is not an
Albanese fibration). Let ${\cal C}=\{\pi_i: S \longrightarrow C_i
\ \mbox{fibrations}, c_i=g(C_i)\geq 2,\pi_i \not= f\}_{i \in I}$.
Assume ${\cal C}\not= \emptyset$ and let $c=\mbox{max} \{c_i \vert
i \in I\}$. Then

(i) $\lambda(f) \geq 4 + \frac{c-1}{g-c}$

(ii) If, moreover, $\mbox{\rm dim}\,alb(S)=1$ (then necessarily $b=0$)
we have
$$\lambda(f)\geq 4 +\frac{q-1}{g-q}.$$
\end{theorem}

\medskip

\vglue1truecm

We want to thank Professor Kazuhiro Konno who kindly communicated
to us the result in Proposition 2.3 (i). The second author would
like to thank also all the staff of the Departament de Matem\`atica
Aplicada I, Universitat Polit\`ecnica de Catalunya for the extremely
warm hospitality received during his visit in February 1998.

\section{Preliminaries}

Here we give a brief run-down of Xiao's method to estimate $\lambda(f)$.
Its method uses a result of Harder and Narasimhan and the
theorem of Clifford.

Let ${\cal{E}}$ be a locally free sheaf on $B$ and let ${\cal{GR}}({\cal{E}})$ be the set of
the locally free subsheaves of ${\cal{E}}$; it is defined a function
$\mu:{\cal{GR}}({\cal{E}})\rightarrow \bbbq$, ${\cal{F}}\mapsto
\mu({\cal{F}})={\rm{deg}}({\cal{F}})/{\rm{rank}}({\cal{F}})$. We recall that ${\cal{E}}$ is
Mumford-stable (respectively Mumford-semistable) if for every proper subbundle
${\cal{F}}$ of ${\cal{E}}$, $0<{\rm{rank}}({\cal{F}})<{\rm{rank}}({\cal{E}})$ we have
$$
\mu({\cal{F}})<\mu({\cal{E}})\,\,\, ({\rm{resp.}}\mu({\cal{F}})\leq\mu({\cal{E}})).
$$
The Harder-Narasimhan theorem concerns the maximum value for $\mu$.

\begin{theorem}\label{hardernarasimhan}
Let ${\cal{E}}$ be a locally free sheaf on a smooth curve $B$, there exists a unique filtration by subbundles
$$
0={\cal{E}}_{0}\subset{\cal{E}}_{1}\subset\cdots\subset{\cal{E}}_{l}={\cal{E}}
$$
\noindent
such that, for $i=1,...,l$, ${\cal{E}}_{i}/{\cal{E}}_{i-1}$
is the maximal semistable subbundle of ${\cal{E}}/{\cal{E}}_{i-1}$.
We put $\mu_{i}=\mu({\cal{E}}_{i}/{\cal{E}}_{i-1})$. In particular
for every $i=1,...,l$, ${\cal{E}}_{i}/{\cal{E}}_{i-1}$ is the
unique subbundle of ${\cal{E}}/{\cal{E}}_{i-1}$
such that for every subbundle ${\cal{F}}$ of
${\cal{E}}/{\cal{E}}_{i-1}$ we have $\mu({\cal{F}})\leq \mu_{i}$
and if
$\mu({\cal{F}})=\mu_{i}$ then
${\cal{F}}\subset{\cal{E}}_{i}/{\cal{E}}_{i-1}$.
Moreover $\mu_{1}>\mu_{2}>...>\mu_{l}.$
\end{theorem}
\proof See \cite{HN}.\qed

The numbers
$\{\mu_{i}\}_{1\leq i\leq l}$ are called
the Harder-Narasimhan slopes of ${\cal{E}}$.

We shall use the following result that relates the Harder-Narasimhan
filtration of a direct sum to the filtrations of the summands.

\begin{proposition}\label{directsum}
Let ${\cal{E}}$, ${\cal{H}}$, ${\cal{K}}$ be locally
free sheaves on a smooth curve $B$. Let
$0={\cal{E}}_{0}\subset{\cal{E}}_{1}\subset\dots\subset{\cal{E}}_{\ell}$,
$0={\cal{H}}_{0}\subset{\cal{H}}_{1}\subset\dots\subset{\cal{H}}_{\ell_{1}}$,
$0={\cal{K}}_{0}\subset{\cal{K}}_{1}\subset\dots\subset{\cal{K}}_{\ell_{2}}$ their
Harder-Narasimhan filtrations and $\{\mu_{i}\}_{1\leq i\leq \ell}$,
$\{\mu^{{\cal{H}}}_{i}\}_{1\leq i\leq \ell_{1}}$,
$\{\mu^{{\cal{K}}}_{i}\}_{1\leq i\leq\ell_{2}}$ their Harder-Narasimhan slopes. Assume
${\cal{H}}\oplus{\cal{K}}={\cal{E}}$. Then we can define
$\psi :\{0,\cdots ,\ell\}\rightarrow \{0,1,\cdots ,\ell_{1}\}$, $\phi :\{0,\cdots
,\ell\}\rightarrow
\{0,1,\cdots ,\ell_{2}\}$ such that

$(i)$ $\psi (0)=\phi(0)=0$; for $1\leq i\leq \ell$, $\psi (i)=\psi(i-1)$ if
$\mu^{\cal{H}}_{t}\neq\mu_{i}$ for every
$t\in \{1,\cdots ,\ell_{1}\}$, (respectively, $\phi(i)=\phi(i-1)$ if
$\mu^{\cal{K}}_{s}\neq\mu_{i}$ for every $s\in \{1,\cdots ,\ell_{2}\}$) and  $\psi (i)=t$ if
$\mu^{\cal{H}}_{t}=\mu_{i}$,
$(respectively, \phi (i)=s$ if $\mu^{\cal{K}}_{s}=\mu_{i}$);

$(ii)$ ${\cal{E}}_{i}={\cal{H}}_{\psi (i)}\oplus{\cal{K}}_{\phi (i)}$.
\end{proposition}
\proof
Call $\pi_{\cal{H}}:{\cal E}\longrightarrow {\cal H}$,
$\pi_{\cal{K}}:{\cal E}\longrightarrow {\cal K}$ the natural projections.
Let ${\cal E}_1^{\cal{H}}=\pi_{\cal{H}}({\cal E}_1)$, ${\cal E}_1^{\cal{K}}=\pi_{\cal{K}}({\cal
E}_1)$; both are locally free since they are torsion free
(${\cal E}_1^{\cal{H}}  \subseteq {\cal{H}}$, ${\cal E}_1^{\cal{K}} \subseteq {\cal K}$).
We have ${\cal E}_1 \subseteq {\cal E}_1^{\cal{H}} \oplus {\cal E}_1^{\cal{K}}$.

Assume ${\cal E}_1^{\cal{H}}\not= 0$. Since ${\cal E}_1$ is semistable and
${\cal{E}}^{\cal{H}}_{1}$ is a quotient, we have that $\mu ({\cal E}^{\cal{H}}_{1}) \geq
\mu({\cal E}_1)=\mu_1$. From the inclusions
${\cal E}_1^{\cal{H}} \subseteq {\cal H} \subseteq {\cal E}$ we get
$\mu_1 \leq \mu({\cal E}_1^{\cal{H}})\leq \mu_1^{\cal{H}}\leq \mu_1$
 since ${\cal H}_1,
{\cal E}_1$ are the maximal unstabilizing sheaves in ${\cal H}$ and
${\cal E}$ respectively. Hence $\mu_1=\mu_1^{\cal{H}}$ and ${\cal E}_1^{\cal{H}}
\subseteq {\cal H}_1 \subseteq {\cal E}_1$
 by the maximality of ${\cal H}_1$ and of ${\cal E}_1$ (see \ref{hardernarasimhan}). The
same argument works if ${\cal E}_1^{\cal{K}} \not= 0$.

Assume $\mu_{1}^{\cal{K}} \not= \mu_1$. Then necessarily ${\cal E}_1^{\cal{K}}
=0$ and ${\cal E}_1^{\cal{H}} \not= 0$. Hence $\mu_1^{\cal{H}}=\mu_1$, ${\cal E}_1
\subseteq {\cal E}_1^1 \subseteq {\cal F}_1$ and then ${\cal E}_1 = {\cal F}_1
={\cal H}_1 \oplus {\cal K}_0$ by maximality. The same argument works if
$\mu_1^{\cal{H}} \not= \mu_1$.

Assume $\mu_1^{\cal{H}} =\mu_1^{\cal{K}} =\mu_1$. Then ${\cal E}_1 \subseteq {\cal
E}_1^{\cal{H}}
\oplus {\cal E}_1^{\cal{K}} \subseteq {\cal H}_1 \oplus {\cal K}_1$ with
$\mu({\cal H}_1 \oplus {\cal K}_1)=\mu_1$. Again by maximality of
${\cal E}_1$ we conclude ${\cal E}_1={\cal H}_1 \oplus {\cal K}_1$.

The proof concludes by induction dealing with
${\cal E}/{\cal E}_1 = {\cal H}/{\cal H}_{\psi(1)} \oplus {\cal K}/
{\cal K}_{\phi(1)}$.\qed\newline
\begin{corollary} \label{directsum2}
With the above notations we have:
$$\begin{array}{l}
\mbox{\rm max}\{\mu_1^{\cal{H}},\mu_1^{\cal{K}}\}=\mu_1\\
\mbox{\rm min}\{\mu^{\cal{H}}_{\ell_1},\mu^{\cal{K}}_{\ell_2}\}=\mu_\ell\\
\end{array}
$$
\end{corollary}
\proof
Obvious.
\qed

\medskip

We will use the well-known refined version of Clifford's theorem:

\medskip

\begin{theorem}{\bf{Clifford-plus.}}\label{clifford}
Let $F$ be a smooth curve of genus $g$. Let $\eta$ be a divisor of degree $d$ such that the linear system $\mid \eta\mid$ has dimension $r-1$ and let
$\phi_{\mid \eta\mid} : F\rightarrow I\!\!P^{r-1}$ be the rational map associated to $D$. Then it holds:

\begin{itemize}
\item[(i)] if $d\leq 2g-2$ then $d\geq 2r-2$;

\item[(ii)] if ${\rm{deg}}(\phi_{\mid \eta\mid})=m$ then $d\geq m(r-1)$;

\item[(iii)] if $\phi_{\mid \eta\mid}=1$ then $a)$ if $d\leq g-1$ then $d\geq 3r-4$, $b)$ if $d\geq g$ then $d\geq \frac{1}{2}(3r+g-4)$.\\
\end{itemize}

Moreover if there exists a smooth curve $C$, a double cover
$\sigma : F\rightarrow C$ and a divisor $\eta$ on $C$ such that
$\mid \eta\mid =\zeta+\sigma^{\star}\mid\eta^{'}\mid$, where $\zeta$ is the fixed part of
$\mid \eta\mid$ then $d\geq 2(r-1+g(C))$.

\end{theorem}
\proof See \cite{Be}[Lemme 5.1]\qed

\medskip

Assume that ${\cal{O}}_{S}(H)$ is an invertible sheaf on $S$ such that
${\cal{E}}=f_{\star}{\cal{O}}_{S}(H)$ is a rank $g$, locally free sheaf on $B$. Let
$0={\cal{E}}_{0}\subset{\cal{E}}_{1}\subset\cdots\subset{\cal{E}}_{l}={\cal{E}}$ be its
Harder-Narasimhan filtration and let $\eta\in{\rm{Pic}}(B)$ be a sufficiently ample
sheaf such that ${\cal{E}}_{i}(\eta)={\cal{E}}_{i}\otimes\eta$ is globally generated.
The natural sheaf homomorphism $f^{\ast}{\cal{E}}_{i}\rightarrow f^{\ast}f_{\ast}
{\cal{O}}_{S}(H)\rightarrow{\cal{O}}_{S}(H)$ induces a rational map
$\rho_{i}:S\rightarrow I\!\!P({\cal{E}}_{i})$. Let $\sigma:\tilde{S}\rightarrow S$ be
the elimination of the indeterminacy of $\rho_{i}$ for every $i$ and let $I\!\!P_{i}$ be the
sublinear system of
$\mid\sigma{\ast}( H+f^{\star}\eta)\mid$ such that
$I\!\!P_{i}=I\!\!P(H^{0}(B,{\cal{E}}_{i}(\eta)))$ for $i=1,...,l$. Xiao defined the following
divisors on $\tilde{S}$:
$Z_{i}$ which is the fixed part of
$I\!\!P_{i}$,
$M_{i}=\sigma{\ast}H - Z_{i}$ which is the moving part of $I\!\!P_{i}$ and
$N_{i}=M_{i}-\mu_{i}F$. By
\cite{Na} for every $i$ $N_{i}$ is a nef $\bbbq$-divisor. We observe that the restriction
$I\!\!P_{i_{|F}}$ is a sublinear system of $\mid H_{\mid F}\mid$ of dimension at least
$r_{i}={\rm{rk}}({\cal{E}}_{i})$, with fixed part $Z_{i_{|F}}$, moving part $M_{i_{|F}}$ and
of degree $d_{i}=N_{i}F$. It is easy to see that these definitions do not depend on
$\eta$. Thus we can give the following definition:

\begin{definition}{\rm{
We call $\{M_{i|F}, r_{i},d_{i}\}$ the Xiao's
data associated to the sheaf ${\cal{E}}=f_{\star}{\cal{O}}_{S}(H)$.}}
\end{definition}

\begin{proposition}\label{xiao}
Let $f:S\rightarrow B$ be a fibration with general fibre $F$. Let $H$ be a divisor on $S$
and suppose there are a sequence of effective divisors on a suitable blow up
$\sigma:\tilde{S}\rightarrow S$, $Z_{1}\geq Z_{2}\cdots\geq Z_{l}\geq Z_{l+1}=0$ and a sequence of rational numbers $\mu_{1}>\mu_{2}>\cdots >\mu_{l} \geq \mu_{l+1}=0$
such that for every $i$, $N_{i}=\sigma^{\ast}H-Z_{i}-\mu_{i}F$ is a nef $\bbbq$-divisor
then
$$
H^{2}\geq\sum_{i=1}^{l}(d_{i}+d_{i+1})(\mu_{i}-\mu_{i+1})
$$
where $d_{i}=N_{i}F$.
\end{proposition}
\proof See \cite{X}[Lemma 2].\qed

\medskip

Xiao's method is to combine \ref{xiao}
and \ref{clifford} when $H=K_{S/B}$ or when $H$
induces on the general fibre $F$ a sublinear system
of $\mid K_{F}\mid$.

\section{The slope of double covers}

\begin{definition}\label{doublecoverfibration}{\rm{
Let $f:S\longrightarrow B$ be a relatively minimal fibration.
We say that $f$ is a {\rm{double cover fibration}} (double cover, for short) if there exists a
relatively minimal fibration $\phi:V \longrightarrow B$ and a rational
map $\pi: S \displaystyle \mathop {--- \rightarrow} V$ over $B$ which is
a generically two to one map. Otherwise we say that $f$ is a {\rm{non double cover fibration}} (non double cover for short).}}
\end{definition}

Roughly speaking, double cover fibrations correspond to the
curves with an involution in the theory of curves.
Of course, if $f$ is a double cover fibration then $F$ is a double
cover. The converse is not true. In \cite{Ba} Example 1.2 a bielliptic
fibration is given for which a non trivial base change is needed in
order to be a double cover fibration. In \cite{BN} is proved that
if the general fibre $F$ is a double cover of a curve of fixed genus
$\gamma$
in a unique way, then the corresponding involution glues to a global
involution of $S$ and so $f$ itself is a double cover. It is easy to
check that if $g \geq 4 \gamma +2$ such condition holds.

Let $f$ be a double cover fibration. We
have
$$\xymatrix{
{\widetilde S} \ar[r]^{\widetilde \pi} \ar[d]^{\sigma}
& {\widetilde V} \ar[d]^{\eta} \\
S \ar[r]^{\pi} \ar[d]^f &V \ar[dl]^{\phi}\\
B
}$$
where $\phi$ and $\pi$ exits by definition, $\pi$ is a generically
2-to-1 rational map and $\phi$ is a relatively minimal fibration;
$\eta:\widetilde V \longrightarrow V$ and $\sigma: \widetilde S
\longrightarrow S$ are any birational maps such that the induced
rational map $\widetilde \pi=\eta^{-1} \circ \pi \circ \sigma$ is a
morphism.  Let
$\widetilde f=f \circ \sigma$ and $\widetilde \phi=\phi \circ \eta$.
Note that at general $t \in B$\ $f^{-1}(t)= \widetilde f^{-1}(t)$,
$\phi^{-1}(t)= \widetilde \phi^{-1} (t)$.
The map $\eta \circ \widetilde \pi$ factorizes by Stein Theorem as
$\eta \circ \widetilde \pi=\pi_0 \circ u$, where $\pi_0$ is finite and
$u$ is birational. Let $R$ be the branching divisor of $\pi_0$ and
${\cal L}\in  Pic(V)$ such that ${\cal L}^{\otimes 2}={\cal O}_V(R)$.
By standard theory of cyclic coverings we have that
$$\begin{array}{rl}
f_{\ast}\omega_{S/B}=\widetilde f_{\ast}\omega_{\widetilde
S/B}=\phi_{\ast}((\eta \circ \widetilde \pi)_{\ast} \omega_{\widetilde
S /B})=\phi_{\ast}(\omega_{V/B} \oplus (\omega_{V/B} \otimes {\cal L}))
=\\
=\phi_{\ast} \omega_{V/B} \oplus \phi_{\ast}(\omega_{V/B}\otimes
{\cal L}).\end{array}$$

\bigskip

\begin{remark}\label{decomposition}{\rm{With the previous notation.
Let ${\cal H}=\phi_{\ast}\omega_{V/B}, {\cal K}=\phi_{\ast} (\omega
_{V/B} \otimes {\cal L})$. Let $h=q(S)-b, s_1=q(V)-b
$ and $s_2=h-s_1$. According to Fujita's decomposition (see\cite{Fu}) we
have $f_{\ast} \omega_{S/B}={\cal E} \oplus {\cal O}^{\oplus s}_{B}$.
A simple computation shows that then we obtain
${\cal H}={\cal F}\oplus {\cal O}^{\oplus s_1}_{B}$,
${\cal K}={\cal G} \oplus {\cal O}^{\oplus s_2}_{B}$ and
${\cal E}={\cal F}\oplus {\cal G}.$}}
\end{remark}

We define $\chi_{1}={\rm{deg}}{\cal{H}}$ and $\chi_{2}=
{\rm{deg}}{\cal{K}}.$

Let $0\subset{\cal{H}}_{1}\subset\dots\subset {\cal{H}}_{{\ell_{1}}-1}
\subset{\cal{H}}_{{\ell_{1}}}$ and $0\subset{\cal{K}}_{1}\subset\dots\subset {\cal{K}}_{l_{2}-1}
\subset{\cal{K}}_{l_{2}}$
be respectively, the Harder Narasimhan filtration of ${\cal{H}}$ and ${\cal{K}}$.
Let $\{(M_{i\mid E}^{{\cal{H}}},r_{i}^{{\cal{H}}}, d_{i}^{{\cal{H}}})\}_{i=1}^{{\ell_{1}}}$ and
$\{(M_{i\mid E}^{{\cal{K}}},r_{i}^{{\cal{K}}}, d_{i}^{{\cal{K}}})\}_{i=1}^{l_{2}}$ be
respectively the Xiao's data of ${\cal{H}}$ and ${\cal{K}}$.

In order to estimate a lower bound for the slope in the double cover
case the following is the key technical result. The first part is due
to Konno (cf. \cite{K2}). We reproduce here a proof for lack of a
suitable reference.

\begin{proposition} \label{konnobound}
Let $f:S\rightarrow B$ be a genus $g$, relatively minimal non isotrivial fibration,
which is a double cover of a genus $\gamma \geq 1$ fibration $\phi :V\rightarrow
B$.
Let $\chi_{f}={\rm{deg}}(f_{\ast}\omega_{S/B})$ and consider the previous
notations.
Assume $g\geq 2\gamma +1$. Then

\begin{itemize}
\item[(i)] $K_{S/B}^{2}-4\chi\geq -4(\mu_{1}^{\cal{H}}+
\mu_{l}^{\cal{H}})+2(g-2\gamma+1){\rm{max}}
\{\frac{\mu_{1}^{\cal{H}}}{\gamma},\mu_{l}^{\cal{H}}\}.$
In particular if $g\geq 4\gamma+1$ then $\lambda(f)\geq 4$.

\item[(ii)]
$K_{S/B}^{2}\geq \frac{8g(g-1)}{g^{2}+g-1}\chi_{1}.$

\item[(iii)] If $g\geq 2\gamma +s_{2}$ then

$$
K_{S/B}^{2}\geq 4\frac{(g-1)(g-s_{2}-1)}{(g-1)(g-\gamma)-s_{2}g}\chi_{2}.
$$

If $g\leq 2\gamma+s_{2}$ then
$$
K_{S/B}^{2}\geq 8\frac{g(g-1)}{g^{2}+g-1)}\chi_{2}.
$$
\end{itemize}
\end{proposition}

\proof We have obtained
$u=\eta\circ{\widetilde\pi}:{\widetilde S}\rightarrow V$ a generically $2$-to-$1$
morphism from a blow-up of $S$ onto a relatively minimal genus $b$
pencil of genus $\gamma$.
Now consider
$$
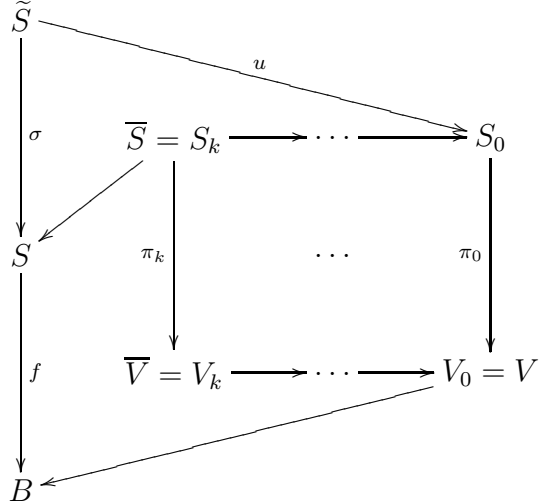
\xymatrix@C=1truecm@R=1truecm{
{\widetilde S}\ar[dd]^{\sigma }\ar[drrr]^{u}&\\
&{\overline S}
=S_{k}\ar[dl]\ar[dd]_{\pi_{k}}\ar[r]&\dots\ar[r]&S_{0}\ar[dd]_{\pi_{0}}\\
S\ar[dd]^{f }&&\dots&\\
&\overline V=V_{k}\ar[r]&\dots\ar[r]&V_{0}=V\ar[dlll]\\
B&&&}
$$

\noindent where:

$\bullet$  $\pi=\pi_{0}\circ u$ is the Stein factorization of $\pi$, with
$u$ birational,  $\pi_{0}$ finite (so it is a double cover) and  $S_{0}$
normal.

 $\bullet$  $\pi_{k}:S_{k}\longrightarrow V_{k}$  is the canonical
resolution of singularities of $\pi_{0}:S_{0}\longrightarrow V_{0}$.

 $\bullet$ $\bar\sigma :S_{k}\longrightarrow S$ is the birational
morphism defined by the relative minimality of $f $. The maps
$\pi_{0}:S_{0}\longrightarrow V_{0}$ and $\pi_{k}:S_{k}\longrightarrow V_{k}$ are
determined by divisors
$R_{0}$  on $V_{0}$,  $R_{k}$  on $V_{k}$ and line bundles ${\cal L}_{0}$, ${\cal L}_{k}$
such that ${\cal L}^{\otimes 2}_{0}={\cal O}_{V_{0}}(R_{0})$,
${\cal L}^{\otimes 2}_{k}={\cal
O}_{V_{k}}(R_{k})$. First of all we have
\begin{equation}\label{eq1:teortresu}
K^{2}_{S/B}-4\chi_f=(K^{2}_{S}-4{\cal X}{\cal O}_{S})-4(b-1)(g-1)\geq
(K^{2}_{\overline S}-4{\cal X}{\cal O}_{\overline S})-4(b-1)(g-1)\, .
\end{equation}

For smooth double covers $\pi_{k}:\overline S\longrightarrow \overline V$ we have
\begin{eqnarray*}
{\cal X}{\cal O}_{\overline S}&=&2{\cal X}{\cal O}_{\overline V}+\frac{1}{2}{\cal L}_{k}K_{\overline
V}+\frac{1}{2}{\cal L}_{k}{\cal L}_{k}\\
K^{2}_{\overline S} &=&2K^{2}_{\overline V}+4{\cal L}_{k}K_{\overline V}+2{\cal L}_{k}{\cal L}_{k}
\end{eqnarray*}
so we have
\begin{equation}\label{eq2:teortresu}
K^{2}_{\overline S}-4{\cal X}{\cal O}_{\overline
S}=2[K^{2}_{V_{k}}-4{\cal X}{\cal O}_{V_{k}}]+2{\cal L}_{k}K_{V_{k}}\, .
\end{equation}
By the canonical resolution of singularities of $\pi_{0}:S_{0}\longrightarrow V_{0}$ we
obtain
\begin{equation}\label{eq3:teortresu}
2[K^{2}_{V_{k}}-4{\cal X}{\cal O}_{V_{k}}]+2{\cal
L}_{k}K_{V_{k}}\geq 2[K^{2}_{V_{0}}-4{\cal X}{\cal O}_{V_{o}}]+2{\cal
L}_{0}K_{V_{0}}
\end{equation}

(i) By (\ref{eq1:teortresu}), (\ref{eq2:teortresu}) and
(\ref{eq3:teortresu}) we have that
$$K_{S/B}^2-4\chi_f \geq 2 (K_{V/B}^2 - 4 \chi_1)+K_{V/B}R$$
where $R=R_0$ is the branch divisor of $S_0 \longrightarrow V_0=V$.

By \ref{xiao} we have a nef $\Rational$-divisor $N_1$ and an
effective divisor $Z_1$ in $V$ such that $K_{V/B}\equiv N_1 + \mu_1^{\cal{H}} E
+Z_1$. Let $R=R_h+R_v$ the decomposition of $R$ in its horizontal and
vertical part respectively. Let $R_h=C_1 + \dots + C_m$ the decomposition
in irreducible components (note that $R$ is reduced since $S_0$ is
normal). Let $n_i$ the multiplicity of $C_i$ in $Z_1$. Then

\begin{eqnarray}
&\sum\limits_{i=1}^m n_i C_i E \leq Z_1 E \leq 2(h-1)
\end{eqnarray}
since $E$ is nef and $Z_1 \leq K_{V/B}$.

Hurwitz formula yields

\begin{eqnarray}
&2(g-2h+1)=R_hE=\sum\limits_{i=1}^mC_iE
\end{eqnarray}

By construction
$$(n_i+1)K_{V/B}-\mu_1^{\cal{H}}E \equiv n_i(K_{V/B}+C_i)+N_1+(Z_1-n_iC_i)$$
we have that

\begin{eqnarray}
&((n_i+1)K_{V/B}-\mu_1E)C_i \geq 0
\end{eqnarray}
since $(K_{V/B}+C_i)C_i \geq 0$ (Hurwitz formula), $N_1C_i \geq 0$
($N_i$ is nef) and  \break
$(Z_1-n_iC_i)C_i\geq 0$ ($C_i$ is not a component of $Z_1-n_iC_i$).

\bigskip

\noindent {\bf Claim.} $K_{V/B}R \geq \frac{2(g-2\gamma+1)}{h} \mu_1^{\cal{H}}$

\noindent Proof of the Claim. We can assume $n_1 \geq n_2 \geq \dots \geq
n_m \geq 0$.

If $h-1 \geq n_1$ ($\geq n_i$ for all $i$) we have that $(hK_{V/B}-
\mu_1E)C_i \geq 0$ by (6) since $K_{V/B}$ is nef.

Assume $h \leq n_1$. Since $n_1 C_1 E \leq 2(h-1)$ we must have
$C_1E=1$. Note that (4) gives $n_i \leq 2h-2-n_1$ for $i \geq 2$. Hence,
using (5) and (6) we have
$$\begin{array}{rl}
K_{V/B}R_h \geq \mu_1^{\cal{H}} \sum\limits_{i=1}^m\frac{1}{n_i+1} C_iE &\geq
\mu_1^{\cal{H}} \left (\frac{C_1E}{n_1+1} + \frac{(R_h-C_1)E}{2h-1-n_1}\right )=\\
&\\
&=\mu{_1}^{\cal{H}} \left (\frac{1}{n_1+1} + \frac{2g-4h+1}{2h-1-n_1}\right )\geq
\mu{_1}^{\cal{H}} \frac{2(g-2h+1)}{h}
\end{array}$$
since $n_1\geq h$.
This proves the Claim.

Finally, since $K_{V/B}-\mu_nE$ is nef we have by (5)
$$K_{V/B}R \geq 2(g-2h+1)\mu_{\ell_1}^{\cal{H}}$$

In \cite{X} p. 460 Xiao gives the following bound for any fibration
$$K_{V/B}^2  \geq 4 \chi_1 -2 (\mu_{1}^{\cal{H}}+\mu_{\ell_1}^{\cal{H}})$$

So
$$K_{V/B}^2-4\chi_f \geq -4(\mu_{1}^{\cal{H}} +
\mu_{\ell_1}^{\cal{H}})+2(g-2\gamma+1)\mbox{max}
\{\frac{\mu_1^{\cal{H}}}{h},\mu_{\ell_1}^{\cal{H}}\}.$$

\medskip

(ii) We consider the Xiao's data
$\{(M_{i\mid E}^{{\cal{H}}},r_{i}^{{\cal{H}}}, d_{i}^{{\cal{H}}})\}_{i=1}^{{\ell_{1}}}$. Since
$\mid M_{i\mid E}^{{\cal{H}}}\mid$ is a sublinear system of $\mid K_{E}\mid$ by Clifford's lemma $d_{i}^{{\cal{H}}}\geq 2(r_{i}^{{\cal{H}}}-1)$. Then it induces on $F$ a linear system of degree $a_{i}\geq 4(r_{i}^{{\cal{H}}}-1)$. Hence, for $1\leq i\leq {
\ell_{1}}-1$ we have $a_{i}+a_{i+1}\geq 8r_{i}^{\cal{H}}-4$. For $i={\ell_{1}}$, $a_{{\ell_{1}}}+a_{{\ell_{1}}+1}\geq 4r_{{\ell_{1}}}-4 +2g-2\geq 8h-4$ since $g\geq 2h+1$ by hypothesis. By \ref{xiao} we obtain
$$
K_{S/B}^{2}\geq 8\sum_{i=1}^{{\ell_{1}}}r_{i}^{{\cal{H}}}(\mu_{i}^{{\cal{H}}}-\mu_{i+1}^{{\cal{H}}})-4\mu_{1}^{{\cal{H}}}=8\chi_{1}-4\mu_{1}^{{\cal{H}}}\geq 8\chi_{1}-4\mu_{1}^{{\cal{E}}}
$$\noindent since $\mu_{1}^{{\cal{E}}}\geq\mu_{1}^{{\cal{H}}}$. Hence
by (i) we have: $(1+\frac{2g-1}{g(g-1)})K_{S/B}^{2}\geq
8\chi_{1}$.

\medskip

(iii) Now we want compare $K_{S/B}^{2}$ with $\chi_{2}$.
Let $N_{i}^{{\cal{K}}}$ and
$Z_{i}^{{\cal{K}}}$ be the divisor on $V$ (on a suitable blow up of $V$) associated to
the Harder-Narasimhan decomposition of ${\cal{K}}$. We put
$N_{i}=\pi^{\ast}(N_{i}^{{\cal{K}}})$,
$Z_{i}=\pi^{\ast}(Z_{i}^{{\cal{K}}})$,
$H=K_{S/B}$, $\mu_{i}=\mu_{i}^
{{\cal{K}}}$ where $i=1,\dots,l_{2}$.
Since $H_{i}=H-Z_{i}-\mu_{i}F=\pi^{\ast}
(K_{V/B}+{\cal{L}}-Z_{i}^{{\cal{K}}}-
\mu_{i}^{{\cal{K}}}E)$ then $H_{i}$ is nef and by
\ref{xiao} $H^{2}\geq\sum_{i=1}^{l_{2}}(d_{i}+d_{i+1})
(\mu_{i}-\mu_{i+1})$ where $d_{i}=N_{i}F=2d_{i}^{{\cal{K}}}$,
$i=1,\cdots,l_{2}$. We consider Xiao'{s} data for ${\cal{K}}$:
$\{(M_{i\mid E}^{{\cal{K}}},r_{i}^{{\cal{K}}}, d_{i}^{{\cal{K}}})\}_
{i=1}^{l_{2}}$. Now the linear systems
$\mid M_{i\mid E}^{{\cal{K}}}\mid$ are sub-linear systems of
$\mid K_{E}+{\cal{L}}_{\mid E}\mid$, so not always they are special.

We put $r_{i}=r_{i}^{{\cal{K}}}$. We have:
$d_{i}^{{\cal{K}}}\geq 2(r_{i}^{{\cal{K}}}-1)$ if $r_{i}^{{\cal{K}}}\leq \gamma$ and
$d_{i}^{{\cal{K}}}=r_{i}^{{\cal{K}}}+\gamma -1$ if $r_{i}^{{\cal{K}}}\geq \gamma-1$. We distinguish two cases: $g\geq 2\gamma+s_{2}$ or $g\leq 2\gamma+s_{2}$.

{\bf{First Case}}:$g\geq 2\gamma+s_{2}$.
If we consider the degree as a function of the rank we easily see that
$d_{i}^{{\cal{K}}}\geq \frac{g-s_{2}-1}{g-s_{2}-\gamma-1}$. Thus
$$
d_{i}^{{\cal{K}}}+d_{i+1}^{{\cal{K}}}\geq 2\frac{g-s_{2}-1}{g-s_{2}-\gamma-1}r_{i}-
\frac{g-s_{2}-1}{g-s_{2}-\gamma-1}
$$
if $i\leq l_{2}-1$ and $
d_{l_{2}}^{{\cal{K}}}+d_{l_{2}+1}^{{\cal{K}}}\geq 2\frac{g-s_{2}-1}{g-s_{2}-\gamma-1}(g-s_{2}-\gamma)-
2\frac{g-s_{2}-1}{g-s_{2}-\gamma-1}.$ We put $A=\frac{g-s_{2}-1}{g-s_{2}-\gamma-1}$ and $B=2\frac{g-s_{2}-1}{g-s_{2}-\gamma-1}$, then

$$
K_{S/B}^{2}\geq\sum_{i=1}^{l_{2}}(d_{i}+d_{i+1})(\mu_{i}-\mu_{i+1})\geq \sum_{i=1}^{l{_2}}(4Ar_{i}-B)(\mu_{i}-\mu_{i+1})-B\mu_{l_{2}},
$$
that is: $K_{S/B}^{2}\geq 4A\chi_{2}-B(\mu_{1}
+\mu_{l_{2}}).$ We recall that $K_{S/B}^{2}\geq
d_{l_{2}}(\mu_{1}+\mu_{l_{2}})$; then
$$
(1+\frac{2g-2-2s_{2}}{(2g-2)(g-\gamma-s_{2}-1)}K_{S/B}^{2}\geq 4A\chi_{2}
$$
so: $K_{S/B}^{2}\geq 4\frac{(g-1)(g-s_{2}-1)}{(g-1)(g-\gamma)-s_{2}g}\chi_{2}.$\newline

{\bf{Second Case}}:$g\leq 2\gamma+s_{2}$.

If $g\leq 2\gamma+s_{2}$ then $d_{i}\geq 4(r_{i}-1)$, so $K_{S/B}^{2}\geq
8\chi_{2} -4\mu_{1}$ then by (i)

$$
K_{S/B}^{2}\geq 8\frac{g(g-1)}{g^{2}+g-1}\chi_{2}.
$$
\qed

\medskip

In the next two theorems we find lower bounds for $\lambda(f)$
in the case of double covers, considering or not the influence
of the relative irregularity of $f$ (see \S 3).
Our limiting functions $l=l(g,\gamma ,s_{1},s_{2})$ or
${\widetilde l}={\widetilde l}(g,\gamma)$
have rather complicate
expressions to be able to check their sharpness. Nevertheless we can give
examples to check its assimptotic good behaviour.

\begin{theorem}\label{boundirregularity} Let $f:S\rightarrow B$ be a genus $g$,
relatively minimal, non isotrivial fibration. Assume $f$ is a double cover fibration of a
fibration of genus $\gamma$. We use the notations of \ref{decomposition}.

\begin{itemize}
\item[(i)] If $g\geq 2\gamma+s_{2}$ and $g>4\gamma+1$ then
$$
\lambda(f)\geq 4+4\frac{(g-4\gamma-1)[(g-1)(\gamma-1)+s_{2}]}{(g-4\gamma -1)[(g-1)(g-\gamma)-gs_{2}]+2(g-1)(g-s_{2}-1)(\gamma -s_{1})\gamma}.
$$
\item[(ii)] If $4\gamma +1 \leq g\leq 2\gamma+s_{2}$ then
$$
\lambda(f)\geq 4 +8\frac{2(g-4\gamma-1)(g-3}{2(g-4\gamma-1)+8(g-1)(\gamma-s_{1})\gamma}.
$$
\end{itemize}
\end{theorem}

\proof
Looking independently to the cases
$\mu_{{\ell_{1}}}^{{\cal{H}}}\geq \mu_{1}^{{\cal{H}}}/\gamma$ and
$\mu_{{\ell_{1}}}^{{\cal{H}}}\leq \mu_{1}^{{\cal{H}}}/\gamma$ in \ref{konnobound} we always get

$$
K_{S/B}^{2}-4\chi\geq \frac{2(g-4\gamma-1)}{\gamma}\mu_{1}^{{\cal{H}}}.
$$
Since $\mu_{1}^{{\cal{H}}}\geq\frac{1}{\gamma-s_{1}}\chi_{1}$ and $g\geq 4\gamma+1$ then
$$ K_{S/B}^2\geq y_{1}=4\chi+ \frac{2(g-4\gamma-1)}{\gamma(\gamma -s_{1})}\chi_{1}.$$
By Proposition 2.3 we know that if $g\geq 2\gamma +s_{2}$ then

$$
K_{S/B}^{2}\geq y_{2}=4\frac{(g-1)(g-s_{2}-1)}{(g-1)(g-\gamma)-s_{2}g}\chi-4\frac{(g-1)(g-s_{2}-1)}{(g-1)(g-\gamma)-s_{2}g}\chi_{1}.
$$
We consider the bounds given above and in Proposition 2.3 as
functions of $\chi_{1}=x$. We consider the region delimited by this three linear inequalities in the plane $(x,y)$.
Since $y_{2}(x_{0})=y_{1}(x_{0})$
implies
$$
x_{0}=\frac{[4g(\gamma-1)+1-\gamma+s_{2}](\gamma-s_{1})\gamma}{2[(g-1)(g-\gamma)-
gs_{2}](g-4\gamma-1)+4(g-1)(g-s_{2}-1)(\gamma-s_{1})\gamma}\chi
$$
then we find
$$
\lambda(f)\geq 4+4\frac{(g-4\gamma-1)[g(\gamma-1)+1-\gamma+s_{2}]}{(g-4\gamma -1)[(g-1)(g-\gamma)-gs_{2}]+2(g-1)(g-s_{2}-1)(\gamma -s_{1})\gamma}.
$$
If $g\leq 2\gamma +s_{2}$ the same argument shows the claim.\qed\newline

\bigskip

\begin{theorem} \label{senzss}Let $f:S\rightarrow B$ be a genus $g$,
relatively minimal, non isotrivial fibration. Assume $f$ is a double cover fibration of a
fibration of genus $\gamma$. Then, if $g\geq 4\gamma+1$ we have
$$
\lambda(f)\geq 4+4\frac{(\gamma-1)(g-4\gamma-1)}{(g-4\gamma -1)(g-\gamma)+2(g-1)\gamma^{2}}.
$$
\end{theorem}
\proof The same argument of \ref{boundirregularity} works.\qed\newline

\bigskip

\begin{corollary}With the above hypotheses, if $q(S)=q(V)=b+\gamma$ and $g\geq
2\gamma+1$ then
$$
\lambda(f)\geq 4\frac{g-1}{g-\gamma}.
$$
\end{corollary}
\proof We have $s_{2}=0$, $s_{1}=s=\gamma$.
Then the claim follows from \ref{boundirregularity}.\qed\newline

\begin{remark} {\rm{We notice that if $\gamma=1$ we find $\lambda(f)\geq 4$ for biellitic
fibrations of genus $g\geq 5$; see \cite{Ba}.}}
\end{remark}

We give now some examples that show that under extra assumptions, the
bounds given are assimptotically sharp.

\bigskip

\begin{example}\label{ana}
{\rm{Let $A$ be an abelian surface with a
base point free linear system $\vert C \vert$, $C^{2}=4$.
Then g(C)=3. Take $C_{1}, C_{2}$ two smooth and transversal members
and let $\sigma: {\widetilde A} \longrightarrow A$ be the blow-up
at the 4 base points. We have then a fibration $\tau: {\widetilde A}
\longrightarrow {\proj}^{1}$ with general fibre ${\widetilde C}$ a
curve of genus 3. Let $E=E_{1}+E_{2}+E_{3}+E_{4}$
 be the $\sigma$-exceptional reduced and irreducible divisor.
Note that $\tau(E)={\proj}^{1}$.
Let $n>>m>>0$ and $\delta=n{\widetilde C}+mE$; then $\vert 2 \delta
\vert$ has no base point.

We can take then a smooth member $R \in \vert 2 \delta \vert$ and
consider the associated double cover $\pi : S \longrightarrow
{\widetilde A}$. Let $f: S \longrightarrow {\proj}^{1}$ be the
induced fibration, and let $F$ be a general fibre. Then $F$ is
a double cover of ${\widetilde C}$ and we have (note that
$K_{{\widetilde A}}=\sigma^{\ast}K_{A}+E=E$)

\begin{eqnarray}
&&h=g({\widetilde C})=3\nonumber\\
&&g=g(F)=4m+5\nonumber\\
&&K^{2}_{S}=2(K_{{\widetilde A}}+\delta)^{2}=8(m+1)(2n-m-1)\nonumber\\
&&\chi {\cal O}_{S}=2\chi {\cal O}_{{\widetilde A}}+\frac {1}{2}
\delta K_{{\widetilde A}}+ \frac {1}{2} \delta ^{2}=2m(2n-m)+
2(n-m)\nonumber
\end{eqnarray}

Moreover, observe that $\delta$ is nef and big since $\delta^{2}>0$
and $\vert 2 \delta \vert$ moves without base points. Then we can
apply Kawamata-Viehweg vanishing theorem and get that
$h^{1}({\widetilde A},{\cal O}_{{\widetilde A}}(- \delta))=0$; hence

$$q(S)=q({\widetilde A})=2$$

Finally we obtain

$$\lambda(f)= \frac {8(m+1)(2n-m-1)+32(m+1)}{2m(2n-m)+2(n-m)+8(m+1)}$$

Fixing $m$ and making $n$ as big as needed we obtain fibrations with
$g=4m+5, h=3$ and the slope arbitrarily near to

$$\frac {16(m+1)}{4m+2}=4+ \frac {8}{g-3}=\lambda^{'}.$$

Now if in Theorem 3.4 we put the data of the above fibration: $s_{2}=0$,
$s_{1}=2$, $\gamma=3$,
$g=4m+5$ with $m\gg 0$ then we obtain
$$\lambda(f)\geq
4+\frac{8}{(g-3)+\frac{6(g-1)}{g-13}}\simeq 4+\frac{8}{g+3}.
$$}}

\end{example}

\begin{remark}{\rm{
Let $Y$ be a smooth surface, let $B$ be a smooth curve and
denote by
$\pi_{Y}:Y\times B\rightarrow Y$, $\pi_{B}:Y\times B\rightarrow B$ the two natural
projections. Let $C,E\in{\rm{Div}}(Y)$ where $E$ has genus $\gamma$ and $\beta,\eta
\in{\rm{Div}}(B)$ such that there exist $W$, $V$  smooth divisors, $W\in\mid
2(\pi_{Y}^{\ast}(C)+\pi_{B}^{\ast}(\beta))\mid$ and
$V\in\mid\pi_{Y}^{\ast}(E)+\pi_{B}^{\ast}(\eta)\mid$. Let $\tau:Z\rightarrow Y\times B$
be the double cover branched on $W$ and $S=\tau^{\ast}(V)$. We put $\pi=\tau_{\mid
S}$, $\phi=\pi_{B\mid V}$, $f=\phi\circ\pi$ and we assume that $W_{\mid V}=\Delta$ is a
smooth divisor. Then $S$ is smooth, $f:S\rightarrow B$ is a fibration and if
$m={\rm{deg}}(\eta)$, $n={\rm{deg}}(\beta)$ and $g$ is the genus of the general fibre $F$
of $f$ then
$$
\lambda(f)=6\frac{4(g-1)(n+m)+2m(K_{Y}+C+E)^{2}}{(3n+6m)(g-1)+
12m\chi({\cal{O}}_{Y})+6m(g(C)-1) +3nEC}
$$

If we consider a $K3$ surface $Y$ we find that for each smooth curve $B$
and for each $\gamma\geq 2$ there exists a double cover fibration $f$ of genus
$g=4\gamma-3$ such that
$\lambda(f)=\frac{(16n+32m)(\gamma-1)}{(3n+5m)(\gamma-1)+4}$.
If $B=I\!\!P^{1}$ we can take $n=m=1$ and obtain a
genus $5$ fibration $f:S\rightarrow I\!\!P^{1}$  with $\lambda(f)=4$ which is a double
cover of a genus $2$ fibration. In particular the slope $4$ can be achieved by double
cover fibrations with $\gamma >1$.}}
\end{remark}

\begin{remark}\label{stablepoint}{\rm{
The constant $4\frac{g-1}{g-\gamma}$ that appears in Corollary 2.6
plays a curious role in the study of double covers: it appears as
a limit bound when adding fibres to the ramification locus. Indeed,
let $F$ and $E$ be respectively, the general fibre of the fibrations
$f:S\rightarrow B$, $\phi:V\rightarrow B$ where $f$ is assumed to be relatively minimal
and non isotrivial. To simplify we assume (although it is not necessary) that $\pi:S\rightarrow V$ is a double
covering branched on a smooth divisor $\Delta\in\mid 2{\cal{L}}\mid$, where
${\cal{L}}\in{\rm{Pic}}(V)$ such that
$f=\pi\circ\phi$. Let $\Delta_{n}\in \mid 2({\cal{L}}+\phi^{\ast}\eta_{n})\mid$ where
$\eta_{n}$ is a divisor on $B$ of degree $n>0$, let $\pi_{n}:S_{n}\rightarrow V$ be the
double covering branched on $\Delta_{n}$ and $f_{n}=\pi_{n}\circ\phi$. Then the
sequence of the slopes $\{\lambda(f_{n})\}_{n\geq 0}$ is monotonous and
$\lim_{n\rightarrow\infty}\lambda(f_{n})=4\frac{g-1}{g-\gamma}$.

In particular, if $\gamma=0$ or $1$ (hyperelliptic or bielliptic
fibrations)
we obtain $4-\frac{4}{g}$ and $4$ recpectively, which are
the exact lower bounds (cf. \cite{X}, \cite{Ba}).
Nevertheless is not true that $\lambda_{\mbox{exp}}=4\frac{g-1}{g-\gamma}$
is in general a lower bound for double cover fibrations, as the following
example shows.}}

\end{remark}

\begin{example}{\rm{
We start as in \cite{Be}[2.6
Example 3]. Let $Y=A\times H$ where $A$ and
$H$ are elliptic curves $\epsilon$ a point of order two on $A$. Let
$X=Y/\langle\sigma\rangle$ where $\sigma$ is an involution defined on $Y$ by $\sigma
(a,h)=(a+\epsilon,-h)$. We denote by $A^{'}$ the quotient of $A$ by the group
$\{0,\epsilon\}$ and by $h_{1},...,h_{4}$ the points of order $2$ on $H$. Let
$p:X\rightarrow A{'}$ and $q:X\rightarrow I\!\!P^{1}=B$ the two natural elliptic
fibrations on $X$. Clearly $K_{X}=p^{\ast}(\eta)$ where $\eta$ is the divisor on
${\rm{Pic}}^{0}(A^{'})$ associated to the \'etale covering $A\rightarrow A^{'}$. Let
$A^{'}_{i}=q^{-1}(h_{i})$ for $i=1,2,3,4$ and $Q\in A^{'}$. Since the divisor
$\delta=p^{\ast}(\eta+dQ)+A^{'}_{1}$ is 2-divisible on $X$ we consider the double
covering $\mu:V\rightarrow X$ associated to $\delta$. We denote
$p\circ\mu=\phi:V\rightarrow B$ and $q\circ\mu=l:X\rightarrow I\!\!P^{1}$. Let
$C=\phi^{\ast}(P^{'})$ and $E=l^{\ast}(P)$ where $P^{'}\in A^{'}$ and
$P\in B$. Let ${\cal{L}}=nE+mC$, $\pi:S_{d,n,m}\rightarrow V$ the double covering associated to ${\cal{L}}$ and $f_{d,n,m}:S_{d,n,m}\rightarrow B$ the induced fibration on $B$ with fibre $F=\pi^{\ast}(E)$. By the standard theory of double covering we have
:
$$
\lambda(f_{d,n,m})=8\frac{2nm+2nd+5m+5d}{4nm +2nd+5m+6d}
$$
and $g=g(F)=4m+4d+1$, $\gamma=g(E)=2d+1$. In particular
$\lambda_{\rm{exp}}=4\frac{g-1}{g-\gamma}=8\frac{m+d}{2m+d}=
\lim_{n\rightarrow\infty}\lambda(f_{d,n,m})$ but
$\lim_{d\rightarrow\infty}\lambda_{\rm{exp}}=8>
8\frac{2n+5}{2n+6}=\lim_{d\rightarrow\infty}\lambda(f_{d,n,m})$.

We remark that to obtain $\lambda(f)<\lambda_{{\rm{exp}}}$ we have $g\sim
2\gamma$, and that $F$ is a double cover of a double cover.}}
\end{example}

\bigskip

Our techniques requires the assumption $g\geq 4\gamma+1$. However the following
theorem shows that  for a double cover fibration $\lambda(f)\geq 4$
holds with a few
exceptions.

\begin{theorem}\label{conjecture}
Let $f:S\rightarrow B$ be a relatively minimal, non isotrivial double
cover fibration of  $\sigma:V\rightarrow B$. Let $F$ and $E$ be the fibres of $f$
and $\sigma$ respectively and let $g=g(F)$, $\gamma=g(E)$. Assume $F$ is not
hyperelliptic or tetragonal, $\gamma\geq 1$ and $g\geq 2\gamma+11$. Then
$\lambda(f)\geq 4$.
\end{theorem}
\proof
By \cite{Ba}(Theorem 2.1) we can assume $F$ is not bielliptic since $2+11=13\geq 6$. We can
also assume $F$ is not trigonal otherwise
$\lambda(f)\geq\frac{14(g-1)}{3g+1}\geq 4$ if $g\geq 9$ using \cite{K1}
(Main theorem).\newline
Consider the Harder Narasimhan filtration of ${\cal E}=f_
{\ast}\omega_{S/B}$:
$0={\cal{E}}_{0}\subset{\cal{E}}_{1}\subset...\subset{\cal{E}}_{l}={\cal{E}}$ with
slopes $\mu_{1}>...>\mu_{l}\geq 0$ and Xiao's data $\{(M_{i\mid F},r_{i},
d_{i})\}_{i=1}^{l}$ .  Note that if $\mid M_{i\mid F}\mid$ induces a map $\phi_{i}$
we have :

If ${\rm{deg}}(\phi_{i})=1$ $d_{i}\geq 3r_{i}-4$   (if $d_{i}\leq g-1$),
$d_{i}\geq\frac{3r_{i}+g-4}{2}$ (otherwise);

If ${\rm{deg}}(\phi_{i})=2$ $d_{i}\geq 2r_{i}+2$    (since $F$ is not hyperelliptic nor
bielliptic)

If ${\rm{deg}}(\phi_{i})=3$ $d_{i}\geq 3r_{i}$   (since $F$ is not trigonal)

If ${\rm{deg}}(\phi_{i})\geq 4$ $d_{i}\geq 4(r_{i}-1)$ .

Observe that, since $M_{i} \leq M_{i+1}$, the map $\phi_{i}$ factorizes through
$\phi_{i+1}$ and then $d_{i+1}|d_{i}$. Note also that $r_{i+1}\geq r_{i}+1$ and
$d_{i+1}\geq d_{i}$. Then we can prove $d_{i}+d_{i+1}\geq 4r_{i}+1$ with a few
exceptions.  Indeed $\mid M_{i}\mid$ does not define any map only
if $(r_{1},d_{1})=(1,0)$. Then $d_{2}\geq 5=4r_{1}+1$ except if $d_{2}=2,3,4$.  All
these possibilities imply $r_{2}=2$ according to the previous inequalities and
hence $F$ would be hyperelliptic, trigonal or tetragonal, all of these being
impossible by hypothesis. From now on we assume $r_{i}\geq 2$.

If ${\rm{deg}}(\phi_{i})\geq 2$ then $d_{i}\geq 2r_{i}$  and hence
$d_{i}+d_{i+1}\geq 2d_{i}+1\geq 4r_{i}+1$, if $d_i<d_{i+1}$; if
$d_i=d_{i+1}$, then $\phi_i=\phi_{i+1}$ and hence $d_i+d_{i+1}
\geq 4r_i+2$.

If ${\rm{deg}}(\phi_{i})=1$ then also
${\rm{deg}}(\phi_{i+1})=1$. If $d_{i}, d_{i+1}\leq g-1$ then $d_{i}+d_{i+1}\geq
3r_{i}-4+3r_{i+1}-4\geq 6r_{i}-5\geq 4r_{i}+1$ since $r_{i}\geq 3$ ($\phi_{i}$ is
birational).

If $d_{i}\leq g-1$,  $d_{i+1}\geq g$ then $d_{i}+d_{i+1}\geq 2d_{i}+1\geq
6r_{i}-7\geq 4r_{i}+1$ except if $r_{i}=3$.
But then $d_{i}+d_{i+1}\geq (3.3-4)
+g=5+g\geq 13\geq 4r_{i}+1$ since $g\geq 11$ by hypothesis.

Finally assume  $d_{i}, d_{i+1}\geq g$, being $\phi_{i}$ and $\phi_{i+1}$ birational
maps. Then
$$
d_{i}+d_{i+1}\geq \frac{3r_{i}+g-4}{2}+\frac{3r_{i+1}+g-4}{2}\geq
3r_{i}+g-4+\frac{3}{2}\geq 4r_{i}+1
$$\noindent
if $r_{i}\leq g-3$ (the case $r_{i}=g-3$ needs a bit care).

Assume $r_{i}=g-2$. If $r_{i+1}=g$ then $d_{i+1}=2g-2$ and we are done. If
$r_{i+1}=g-1$ then the only case to check is $d_{i}=2g-5$, $d_{i+1}=2g-3$. Note
that then $h^{0}(F,K_{F}-M_{i_{|F}})=h^{0}(F,K_{F}-M_{i+1_{|F}})=1$ since $F$ is not
hyperelliptic. By Riemann-Roch $r_{i}=h^{0}(F,M_{i_{|F}})=1+d_{i}+1-g=g-3$ which
is impossible.

Assume $r_{i}=g-1$. Then $d_{i}=2g-3$ ,
$(M_{i+1_|{F}},d_{i+1})=(r_{l},d_{l})=(g,2g-2)$ and $d_{i}+d_{i+1}=4g-5=4r_{i}-1$.

For $r_{i}=g=r_{l}$ we have $d_{l}+d_{l+1}=2d_{l}=4g-4=4r_{l}-4$. By \ref{xiao} we
conclude
\begin{equation}
\begin{array}{clcl}
K_{S/B}^{2} & \geq & \sum_{i=1}^{l}(d_{i}+d_{i+1})
(\mu_{i}-\mu_{i+1})\geq
& \\

& \geq &
\sum_{i=1}^{l}(4r_{i}+1)(\mu_{i}-\mu_{i+1})-
2(\mu_{l-1}-\mu_{l})-5\mu_{l} & \\
& = & 4\chi+\mu_{1}-2\mu_{l-1}-3\mu_{l}
\end{array}
\end{equation}
\noindent
if $r_{l-1}=g-1$, $d_{l-1}=2g-3$; otherwise
$$
K_{S/B}^{2} \geq
 \sum_{i=1}{l}(4r_{i}+1)(\mu_{i}-\mu_{i+1})-5\mu_{l}
=4\chi+\mu_{1}-5\mu_{l}.
$$

Let us consider first the general case. If  $\mu_{1}\geq 5\mu_{l}$ we are done.
Assume $\mu_{1}\leq 5\mu_{l}$.  Let ${\cal{H}}=\phi_{\ast}(\omega_{V\mid B})$
and ${\cal{K}}=\phi_{\ast}(\omega_{V\mid B}\otimes{\cal{L}})$. By \ref{directsum},
${\cal{E}}_{i}={\cal{H}}_{\psi(i)}\oplus{\cal{K}}_{\phi(i)}$ where
$\mu_{i}=\mu^{{\cal{H}}}_{\psi(i)}=\mu^{{\cal{K}}}_{\phi(i)}$.
By \ref{directsum2}
we have that
$\mu_{l}=
{\rm{min}}\{\mu^{{\cal{H}}}_{{\ell_{1}}},\mu^{{\cal{K}}}_{l_{2}}\}
\leq\mu^{{\cal{H}}}_{{\ell_{1}}}$,
$\mu_{1}={\rm{max}}\{\mu^{{\cal{H}}}_{1},\mu^{{\cal{K}}}_{1}\}
\geq \mu^{{\cal{H}}}_{1}$. Hence  $\mu^{{\cal{H}}}_{1}\leq\mu_{1}<5\mu_{l}\leq
5\mu^{{\cal{H}}}_{{\ell_{1}}}$.  By \ref{konnobound} we have
$$
K_{S/B}^{2}\geq 4\chi -24\mu_{{\ell_{1}}}^{\cal{H}}+2(g-2\gamma+1)
{\rm{max}}
\{\frac{\mu_{1}^{\cal{H}}}{\gamma},\mu_{{\ell_{1}}}^{\cal{H}}\}.
$$
If $\gamma\geq 5$ or
$\frac{\mu_{1}^{\cal{H}}}{\gamma}\leq\mu_{{\ell_{1}}}^{\cal{H}}$ we have ${\rm{max}}
\{\frac{\mu_{1}^{\cal{H}}}{\gamma},\mu_{{\ell_{1}}}^{\cal{H}}\}=\mu_{{\ell_{1}}}^{\cal{H}}$
and hence
$$
K_{S/B}^{2}\geq 4\chi +2(g-2\gamma-11)\mu_{{\ell_{1}}}^{\cal{H}}\geq 4\chi
$$\noindent
when $g\geq 2\gamma+11$.

If $\gamma=2,3,4$ and
$\frac{\mu_{1}^{\cal{H}}}{\gamma}\geq\mu_{{\ell_{1}}}^{\cal{H}}$  then
$$
K_{S/B}^{2}\geq 4\chi
-24\mu_{{\ell_{1}}}^{\cal{H}}+2(g-2\gamma+1)\frac{\mu_{1}^{\cal{H}}}{\gamma}\geq
4\chi +2(g-2\gamma-11)\frac{\mu_{1}^{\cal{H}}}{\gamma}\geq 4\chi
$$\noindent
when $g\geq 2\gamma+11$.

Consider finally the special case $r_{l-1}=g-1$, $d_{l-2}=2g-3$. By \ref{xiao}, with
the notations of \ref{xiao} where $H=\omega_{S/B}$ it follows that $Z_{l-1}$
is a section of $f$ such that $Z_{l-1_{\mid F}}\equiv K_{F}-M_{l-1_{|F}}$.  We recall
that $M_{l-1_{|F}}$ is the base point free linear system induced on the general
fibre by the piece ${\cal{E}}_{l-1}$.  By \ref{directsum} ${\cal{E}}_{l-1}=
{\cal{H}}_{\psi(l-1)}\oplus{\cal{K}}_{\phi(l-1)}$. Since $r_{l-1}=r_{l}-1$ we only
have two possibilities: either ${\cal{H}}_{\psi(l-1)}={\cal{H}}_{{\ell_{1}}}$,
${\cal{K}}_{\phi(l-1)}={\cal{K}}_{l_{2}-1}$ and $r_{l_{2}-1}^{{\cal{K}}}=g-\gamma-1$
or
${\cal{H}}_{\psi(l-1)}={\cal{H}}_{{\ell_{1}}-1}$,
${\cal{K}}_{\phi(l-1)}={\cal{K}}_{l_{2}}$ and $r_{{\ell_{1}}-1}^{{\cal{K}}}=\gamma-1$. We
claim the second possibility can not occur. Indeed consider the double cover
$\pi_{\mid F}:F\rightarrow E$.  We have that
$$
H^{0}(F,\omega_{F})\simeq H^{0}(E,\omega_{E})\oplus
H^{0}(E,\omega_{E}\otimes{\cal{L}}_{\mid E})
$$\noindent
This decomposition means that if $D$ is the ramification divisor on $F$ and
$t\in H^{0}(F,{\cal{O}}_{F}(D))$ then for every $\omega\in
H^{0}(F,\omega_{F})$, $\omega=
t\pi_{\mid F}^{\ast}(\omega_{1}) +\pi_{\mid F}^{\ast}(\omega_{2}) $ where
$\omega_{1}\in  H^{0}(E,\omega_{E})$ and
$\omega_{2}\in H^{0}(E,\omega_{E}\otimes{\cal{L}}_{\mid E})$.

We have $V\in H^{0}(F,\omega_{F})$ a codimension one subspace which produces, after
taking out the base point, the linear series $\mid M_{l-1_{|F}}\mid$. The second
possibility asserts that $V=\pi_{\mid F}^{\ast}V_{1}\oplus
\pi_{\mid F}^{\ast}V_{2}$ where
$V_{2}=H^{0}(E,\omega_{E}\otimes{\cal{L}}_{\mid E})$ and $V_{1}$ is a codimension one
subspace of $H^{0}(E,\omega_{E})$. Since ${\rm{deg}}(\omega_{E}\otimes{\cal{L}}_{\mid E})\geq
2\gamma+10$, $V_{2}$ is base point free. Hence $\pi_{\mid F}^{\ast}V_{2}$ is base
point free: a contradiction since $V$ has a base point.

So we have the following decompositions
$$
{\cal{E}}_{l}={\cal{H}}_{{\ell_{1}}}\oplus{\cal{K}}_{l_{2}},\,\,
{\cal{E}}_{l-1}={\cal{H}}_{{\ell_{1}}}\oplus{\cal{K}}_{l_{2}-1}
$$\noindent
where $r_{l_{2}-1}^{{\cal{K}}}=g-\gamma-1$. If
${\cal{E}}_{l-2}={\cal{H}}_{j}\oplus{\cal{K}}_{k}$ we have several possibilities
according to \ref{directsum}.

If $j={\ell_{1}}$, $k=l_{2}-2$ then
$\mu_{l-1}=\mu({\cal{E}}_{l-1}/{\cal{E}}_{l-2})=\mu
({\cal{K}}_{l_{2}-1}/{\cal{K}}_{l_{2}-2})=\mu_{l_{2}-1}^{\cal{K}}$ and
$\mu_{{\ell_{1}}}^{\cal{H}}>\mu_{l-1}$.

If $j={\ell_{1}}-1$, $k=l_{2}-1$ then
$\mu_{l-1}=\mu
({\cal{H}}_{{\ell_{1}}}/{\cal{H}}_{{\ell_{1}}-1})=\mu_{{\ell_{1}}}^{\cal{H}}$.

If $j={\ell_{1}}-1$, $k=l_{2}-2$ then
$\mu_{l-1}=\mu_{l_{2}}^{\cal{H}}=\mu_{l_{2}-1}^{\cal{K}}$.

In any case we get $\mu_{l-1}=\mu_{{\ell_{1}}}^{\cal{H}}$. Since always happens that
$\mu_{1}\geq\mu_{1}^{\cal{H}}$ and $\mu_{l}\leq\mu_{{\ell_{1}}}^{\cal{H}}$,
(7) reads:
$$
K_{S/B}^{2}\geq 4\chi+\mu_{1}-2\mu_{l-1}-3\mu_{l}\geq
4\chi+\mu_{1}^{\cal{H}}-5\mu_{{\ell_{1}}}^{\cal{H}}.
$$
If $\mu_{1}^{\cal{H}}\geq 5\mu_{{\ell_{1}}}^{\cal{H}}$ we are done.
If $\mu_{1}^{\cal{H}}<5\mu_{{\ell_{1}}}^{\cal{H}}$ then we can repeat the argument of
the general case.
\qed

\section{The slope of non-Albanese fibrations}

In this section we consider the problem of the influence of the
relative irregularity $h=q(S)-b$ on the lower bound of the slope.
In case $f$ is a double cover fibration this problem has been
considered in the previous section so we will deal only with
non double cover fibrations.
The nice fact is that we find a lower bound which is an increasing function of
the genus $g$ and of $h$.

When $h=0$ then the general bound $\lambda(f)\geq 4- \frac{4}{g}$ holds
and is sharp. So we will consider fibrations with $h=q(S)-b>0$. Those
are precisely the fibrations for which the Albanese map of $S$ does not
factorize through $f$ (i.e., $b=0$ and $q>0$ or $S$ is of Albanese
general type). We call such fibrations {\it non-Albanese fibrations}.

Let us first recall the two basic known results in this area:

\begin{theorem}\label{conocido}
Let $f:S\rightarrow B$ be a relatively minimal, non locally trivial,
genus $b$ pencil of curves of
genus $g$.

\begin{itemize}
\item[(i)] If $h=q(S)-b>0$, then $\lambda \geq 4$
\item[(ii)] Let $\mu_{1}={\rm{deg}}{\cal{E}}_{1}/{\rm{rank}}({\cal{E}}_{1})$ where ${\cal{E}}_{1}$ is the maximal semistable subbundles of
$f_{\ast}\omega_{S/B}$. If $g\geq 2$ and $g>q-b$ then $K^{2}_{S/B}\geq\frac{4g(g-1)}{(2g-1)}\mu_{1}$. In particular $\lambda(f) \geq \frac{4g(g-1)}{(2g-1)(g-h)}.$

\end{itemize}
\end{theorem}
\proof $(i)$ is \cite{X}[Theorem 2.4]. $(ii)$ is \cite{K1}[Lemma 2.7]\qed\newline

\begin{theorem}\label{general}
Let $f:S \longrightarrow B$ be a relatively minimal fibration which is not
a double cover fibration. Assume $g=g(F)\geq 5$ and that $f$ is not
locally trivial. Let $h=q(S)-b\geq 1$.

Then

(i) If $h \geq 2$ and $g \geq \frac{3}{2}h+2$ then
$$\begin{array}{rl}
\quad \lambda(f)&\geq \displaystyle{
\frac{8g(g-1)(4g-3h-10)}{8g(g-1)(g-h-2)+3(h-2)(2g-1)}}
\quad \mbox{if}\ F \ \mbox{is not trigonal}\\
&\\
\quad \lambda(f)&\geq \displaystyle{
\frac{4g(g-1)(4g-3h-10)}{4g(g-1)(g-h-2)+(g-4)(2g-1)}}
\quad \mbox{if}\ F \ \mbox{is trigonal}
\end{array}$$

(ii) If $g < \frac{3}{2}h+2$ then
$$\lambda(f)\geq \frac{4g(g-1)(2g-7)}{\frac{4}{3}g(g-1)(g-3)+(g-4)(2g-1)}.
$$
\end{theorem}

\proof

(i), (ii) Consider Fujita's decomposition $f_{\ast}\omega_
{S/B}={\cal A}\oplus {\cal Z}$ with ${\cal Z}={\cal O}_B^{\oplus h}$.
Consider the Harder-Narasimhan filtration of ${\cal A}$:
$$0 = {\cal A}_0 \subseteq {\cal A}_1 \subseteq \dots \subseteq
{\cal A}_\ell =\cal A$$

As in \S 1 we produce nef $\Rational$-divisors $N_i$, and effective
divisors $Z_i$ in a suitable blow-up of $S \quad \sigma:\widetilde S
\longrightarrow S$ such that
$$N_i + \mu_i F +Z_i \equiv N_j + \mu_j F +Z_j \equiv
\sigma^{\ast} K_{S/B}$$
where $\{\mu_i\}$ are the Harder-Narasimhan slopes of ${\cal A}$. Note
that we can define $N_{\ell +1} = \sigma^{\ast} K_{S/B}$, $Z_{\ell +1}=0$,
$\mu_{\ell +1}=0$. Observe also that, if $r_i=\mbox{rk}{\cal A}_i$,
\break
$\sum\limits_{i=1}^\ell r_i(\mu_i-\mu_{i+1})=\mbox{deg}{\cal A}=
\chi_f$.

Each $N_i$ induces on $F$ a base point free linear system of degree $d_i$
and (projective) dimension greater or equal than $r_i -1$. Note that
$N_i + \mu_i F = H_i$ is induced by a map $\varphi_i:S \longrightarrow
{\proj}_B ({\cal A}_i)$ which restricted to fibres induces the above
linear system. By hypothesis $\varphi_i$ is never a double cover onto
the image and so the induced map $\psi_i$ on $F$ is not a double cover.
Hence we have
$$\begin{array}{rl}
d_i &\geq 3(r_i-1) \qquad \mbox{if deg}\psi_i\geq 3\\
&\\
d_i &\geq 3r_i-4 \qquad \mbox{if deg}\psi_i=1 \ \mbox{and} \ d_i\leq g-1\\
&\\
d_i &\geq \frac{3r_i+g-4}{2} \qquad \mbox{if deg}\psi_i=1\ \mbox{and}\
d_i \geq g
\end{array}$$
the latest two inequalities being ``Clifford plus" Lemma.
Considering the above inequalities in the $(r,d)$-plane, we
have the following two possibilities (note that the lines $d=3r-4$ and
$d=\frac{3r+g-4}{2}$ meet exactly at the point ($r=\frac{1}{3}(g+4),
d=g$)) depending on rank${\cal A}=g-h$.

\noindent {\bf Case 1.-} $g-h \geq \frac{1}{3}(g+4)$

In this case note that for every $1 \leq i \leq \ell$, $d_i \geq
\frac{2g-\frac{3}{2}h-5}{g-h-2}r-\frac{g-4}{g-h-2}$ (this border line
joining the point (2,3) and the point $(g-h,2g-\frac{3}{2}h-2)$) except
if $(r_1,d_1)=(1,0)$. Note that $g-h-2>0$ since $g \geq  \frac{3}{2}h+2$.

Note also that by definition we have $d_{\ell +1}=2g-2$. So for
$ 1 \leq i \leq \ell$ we get (since $r_{i+1}\geq r_i+1$)
$$d_i+d_{i+1}\geq \frac{4g-3h-10}{g-h-2}r_i - \frac{3(h-2)}{2(g-h-2)}
=:Ar_i+B$$
except if $(r_1,d_1)=(1,0)$ and $(r_2,d_2)=(2,3)$. In this exceptional
case we get
$$d_1+d_2-Ar_1-B=3-A-B=-\frac{g-\frac{3}{2}h-1}{g-h-2}$$

If this happens $F$ is trigonal since has a linear system of
degree 3 and dimension 1.

Applying Xiao's formula we get, in the general case,
$$\begin{array}{rl}
K_{S/B}^2 &\geq \sum\limits_{i=1}^\ell (d_i+d_{i+1})
(\mu_i-\mu_{i+1})\geq \sum\limits_{i=1}^\ell A r_i (\mu_i - \mu_{i+1})
+\sum\limits_{i=1}^\ell B(\mu_i - \mu_{i+1})=\\
&\\
&=A \chi_f + B \mu_1=\frac{4g-3h-10}{g-h-2} \chi_f -
\frac{3(h-2)}{2(g-h-2)}\mu_1\end{array}$$

Applying again Konno's bound:
$$K_{S/B}^2\geq \frac{4g(g-1)}{2g-1}\mu_1$$
we can eliminate $\mu_1$ and get
$$K_{S/B}^2\geq \frac{8g(g-1)(4g-3h-10)}{8g(g-1)(g-h-2)+3(h-2)(2g-1)}
\chi_f$$

Note that this bound is a strictly increasing function of $h$ and that
$K_{S/B}^2\geq 4\chi_f$ if $h \geq 2$.

\smallskip

In the exceptional case (when $F$ is trigonal) we get
$$K_{S/B}^2 \geq A \chi_f + B \mu_1 - \frac{g-\frac{3}{2}h-1}{g-h-2}
(\mu_1 - \mu_2)\geq A \chi_f + \left (B - \frac{g-\frac{3}{2} h-1}{g-h-2}
\right )\mu_1$$

The same argument using $K_{S/B}^2 \geq \frac{4g(g-1)}{2g-1}\mu_1$ yields
$$K_{S/B}^2 \geq \frac{4g(g-1)(4g-3h-10)}{4g(g-1)(g-h-2)+(g-4)(2g-1)}
\chi_f$$
which is also a strictly increasing function of $h$. In this case we
need $h \geq 4$ to get $K_{S/B}^2 \geq 4 \chi_f$.

\noindent {\bf Case 2.-} $g-h \leq \frac{1}{3} (g+4)$

Let $\overline h=\left [\frac{2}{3} g - \frac{4}{3}\right ]$. Under
our hypotheses $h \geq \overline h$, so we can take $\overline {\cal A}
={\cal A}\oplus  {\cal O}_B^{\oplus (h - \overline h)}$
instead of ${\cal A}$. Hence we get according to whether we are in the general
or in the special case
$$\begin{array}{rll}
K_{S/B}^2&\geq \frac{8g(g-1)(4g-3\overline h -10)}{8g(g-1)(g-\overline h
-2)+3(\overline h -2)(2g-1)} \chi_f &\geq
\frac{8g(g-1)(2g-7)}{\frac{8}{3}g(g-1)(g-3)+(2g-9)(2g-1)}\chi_f\\
&\\
K_{S/B}^2&\geq \frac{4g(g-1)(4g-3\overline h -10)}{4g(g-1)(g-\overline h
-2)+(g-4)(2g-1)} \chi_f &\geq
\frac{4g(g-1)(2g-7)}{\frac{4}{3}g(g-1)(g-3)+(g-4)(2g-1)}\chi_f
\end{array}$$
since both expressions are increasing functions of $h$ and $\overline h
\geq \frac{2}{3} g-1$. Note that the second bound is slightly smaller than
the first one.
\qed

\bigskip

\begin{remark}
{\rm{In the case (iii) of the theorem we could consider that for $1 \leq i
\leq \ell$, $d_i \geq 3r_i-4$ and hence $d_i+d_{i+1} \geq 6r_i-5$
for $1 \leq i \leq \ell-1$. But for $i = \ell$ we would have $d_\ell
+d_{\ell  +1}\geq 2 d_\ell +1 \geq 6 r_\ell -7$ which produces
$$K_{S/B}^2 \geq 6 \chi_f - (5\mu_1+2\mu_\ell)$$

Hence using Xiao's inequality with indexes $\{1,\ell\}$
$$K_{S/B}^2  \geq (d_1 + d_\ell)(\mu_1 - \mu_\ell)+(d_\ell + d_{\ell +1})
\mu_\ell \geq d_\ell(\mu_1 + \mu_\ell)\geq (3g-3h-4)(\mu_1+\mu_\ell)$$
we get
$$K_{S/B}^2 \geq 6 \frac{3g-3h-4}{3g-3h+1}\chi_f$$
which depends on $h$ and is better than (iii) for some special values of
$(g,h)$ but is a {\it decreasing} function of $h$.

Nevertheless we must have in mind that case (iii) of the theorem
is doubtfull to happen. Indeed, by a conjecture of Xiao (cf. \cite{X2})
the following inequality should hold: $h=q_f \leq \frac{1}{2}(g+1)$. This
inequality is true when $b=0$ but is known to be false in general (cf.
\cite{Pi}) although it seems that only the constant term should be
modified.}}

\end{remark}

\bigskip

\begin{remark}{\rm{
In the above theorem we worked with ${\cal Z}={\cal O}_B^{\oplus (q(S)-b)}
$ and $h=\mbox{rank}{\cal Z}$. In most
parts of the proof we only use that deg${\cal
Z}=0$. Hence, we get the same bounds in (ii) if we define $h$ to be the
rank of the degree zero part in Fujita's decomposition of ${\cal E}=
f_{\ast}\omega_{S/B}$ ($h \geq q(S)-b$). Note that then the argument
of Theorem 3.2 (ii) does not work since we do not
know whether ${\cal Z}$
can be cut in pieces of the length we need.
In any case the bound of the previous remark holds for this new
definition of $h$.}}

\end{remark}

\bigskip

\begin{remark}{\rm{
Remember that if $F$ is trigonal we have (cf. \cite{K1} and \cite{SF}):
$$\lambda(f)\geq \frac{14(g-1)}{3g+1}$$
which is better that Theorem 4.19 (ii) (special case) for $g\gg h=q-b$
and that gives $\lambda(f) \geq 4$ if $g \geq 9$.}}
\end{remark}

\bigskip

\begin{remark}{\rm{
As a function on $g$ (fixing $h$) the bounds of \ref{general}
tend to be 4
when $g$ grows (compare Theorem 3.1 (ii) where this limit is 2).}}
\end{remark}

\begin{example}{\rm{
Let $Y$ be a smooth surface,
let $B$ be a smooth curve of genus $b$, $Z=Y\times B$ and let $\pi_{Y}:Z\rightarrow Y$,
$\pi_{B}:Z\rightarrow B$ be the natural projections. If $F\in{\rm{Div}}
(Y)$ is smooth of
genus $g$, $\eta\in{\rm{Div}}(B)$ and there exists an ample and smooth divisor
$S\in\mid \pi_{Y}^{\ast}(F)+\pi_{B}^{\ast}(\eta)\mid$ then the slope of the fibration
$f:S\rightarrow B$ induced on $S$ by $\pi_{B}$ is
$$
\lambda(f)=\frac{6g-6+K_{Y}^{2}+K_{Y}F}{\chi({\cal{O}}_{Y})+g-1}.
$$

Now if $\rho:Y=I\!\!P({\cal{E}})\rightarrow C$ is a ruled
surface, $H$ is a
section such that $H^{2}={\rm{\deg}}({\cal{E}})$ and $F\equiv 3H$
the fibration $f_{m}:S_{m}\rightarrow B$ ($m= \mbox{deg}\eta$) has slope:
$$
\lambda(f_{m})=\frac{15m+16(g(C)-1))}{3m+2(g(C)-1))}
$$
and verifies that $h=q(S_m)-b=g(C)$.
In particular
$\lambda(f_{m})\geq 5$ and $\lim_{m\rightarrow\infty}\lambda(f_{m})=5$.
This result indicates that for any $h$, a general lower bound of
$\lambda(f)$ is below 5.}}

\end{example}

From Theorem \ref{general} we obtain that $\lambda(f)$ controls the existence of other fibrations on $S$:

\bigskip

\begin{theorem}
Let $f:S \longrightarrow B$ be a relatively minimal, non locally trivial
fibration. Let $F$ be a fibre of $f$, $g=g(F)$ and $q=q(S)$.
Assume $f$ is not a
double cover fibration and that $h=q-b\geq1$ (i.e., $f$ is not an
Albanese fibration). Let ${\cal C}=\{\pi_i: S \longrightarrow C_i
\ \mbox{fibrations}, c_i=g(C_i)\geq 2,\pi_i \not= f\}_{i \in I}$.
Assume ${\cal C}\not= \emptyset$ and let $c=\mbox{max} \{c_i \vert
i \in I\}$. Then

(i) $\lambda(f) \geq 4 + \frac{c-1}{g-c}$

(ii) If, moreover, $\mbox{\rm dim}\,alb(S)=1$ (then necessarily $b=0$)
we have
$$\lambda(f)\geq 4 +\frac{q-1}{g-q}$$
\end{theorem}

\proof
Remember that if $f$ is not an Albanese fibration then either
$\mbox{dim}\,alb (S)=2$ or $b=0$ (provided $q(S)\not=0$).

Let $\pi:S \longrightarrow C$ be the fibration with maximal base genus
$c \geq 2$ (if $\mbox{dim}\,alb(S)=1$, then $c=q$ and $\pi=alb$).

Since in any case $f^{\ast} Pic^0(B)$ does not include $\pi^{\ast}
 Pic^0(C)$  we can choose for $n\gg0$, a $n$-torsion element
${\cal L}\in  Pic^0(C)$ such that $\pi^{\ast}{\cal L}^{\otimes i}
\notin f^{\ast} Pic^0(B)$ for $1 \leq i \leq n-1$. Consider the
base change
$$\xymatrix{
\widetilde S \ar[r]^{\widetilde \pi} \ar[d]^{\widetilde \alpha}&
\widetilde C \ar[d]^\alpha \\
S \ar[r]^\pi \ar[d]^f & C\\
B
}$$
and let $\widetilde f=f \circ \widetilde \alpha$. Since ${\cal L}^{\otimes
i}_{\vert F} \not= {\cal O}_F$ for $ 1 \leq i \leq n-1$, $\widetilde f$
has connected fibres and so $\widetilde f$ is again a fibration over $B$.
Let $\widetilde F$ be the fibre of $\widetilde f$. Then if $\widetilde g=
g(\widetilde F)$,
$$\widetilde g-1=n(g-1)$$

Moreover we have
$$q(\widetilde S)=h^1(\widetilde S,{\cal O}_{\widetilde S})=
h^1(S,{\cal O}_S)+ \sum\limits_{i=1}^{n-1} h^1(S,(\pi^{\ast}{\cal L}^{-i})
)$$

>From the exact sequence
$$0 \longrightarrow H^1(B,{\cal L}^{-i}) \longrightarrow H^1
(S,\pi^{\ast}{\cal L}^{-i})\longrightarrow H^0(B,(R^1\pi_{\ast}{\cal O}_S)
\otimes {\cal L}^{-i})\longrightarrow 0$$
and using that $h^0(B,(R^1\pi_{\ast}{\cal O}_S)\otimes {\cal L}^{-i})=0$ except
for a finite number of sheaves ${\cal L}^{-i}\in  Pic^0(C)$ (which can be
avoided with the election of ${\cal L}$ (see \cite{Beau}[Lemme 3.1]
and \cite{BT} \S 3))
we get
$$\widetilde h= q (\widetilde S)-b=q(S)-b+(n-1)(c-1)=h+(n-1)(c-1)$$
since $h^1(B,{\cal L}^{-i})=c-1$ by Riemann-Roch. In particular,
$\widetilde h\geq 2$ if $n\geq 2$.

It is easy to check that if $F$ is trigonal then $\widetilde F$ is
not if $n\gg0$ (see
for example \cite{BT}, Lemma 5.12). On the other hand
$$\lim_{n \rightarrow \infty} \frac{\widetilde g}{\widetilde h}=
\frac{g-1}{c-1}\geq 2$$
since the map $\pi_{\vert F}:F \longrightarrow C$ is at least of degree
two (if it were of degree 1 clearly $F \cong C$ and $S = B \times C$).
Hence if $n\gg0$ the case $\widetilde g< \frac{3}{2} \widetilde h+2$ can
not occur.

So if $n \gg0$ we are under the hypotheses of Theorem 3.2 (ii) (non
trigonal case). Using that the slope is invariant under \'etale changes
of $S$ (cf. \cite{X}) we get
$$\lambda(f)=\lambda(\widetilde f)\geq \frac
{8 \widetilde g(\widetilde g-1)(4 \widetilde g-3\widetilde h-14)}
{8 \widetilde g(\widetilde g-1)(\widetilde g-\widetilde h -3)+5
(\widetilde h-2)(2\widetilde g-1)}$$
for $\widetilde g=n(g-1)+1$, $\widetilde h=h+(n-1)(c-1)$ and $n \in
\Natural$ arbitrarily large. So we can take limit as $n$ grows and get
$$\lambda(f)\geq 4 + \frac{c-1}{g-c}$$

In case dim\,$alb(S)=1$ then clearly $c=q$. Note that if this happens
and $b \geq 1$, then $alb(S)=B$ by the universal property of
Albanese variety.
\qed

\bigskip

\begin{corollary}\label{adria1}
Let $f:S \longrightarrow B$ be as in Theorem 4.19. Assume $\lambda(f)<
4 + \frac{1}{g-2}$. Then $S$ has no other fibration onto a curve of
genus greater or equal than two.
\end{corollary}

\bigskip

\begin{corollary}\label{adria2}
Let $S$ be a minimal surface with $q(S)\geq 2$ and $F \subseteq S$
an irreducible curve of geometric genus $g$. Assume $h^0(S,{\cal O}_S
(F))\geq 2$ and let $f:\widetilde S \longrightarrow {\proj}^1$ be a
relatively minimal fibration with fibre $F$. If $F$ is not a double
cover and $\lambda(f) < 4 + \frac{q-1}{g-q}$ then $S$ is of Albanese
general type.
\end{corollary}

\end{document}